\documentclass[12pt,oneside,reqno]{amsart}
\usepackage{amsmath,amsfonts,amssymb}
\usepackage{cite}

\theoremstyle{plain}

\theoremstyle{definition}

\theoremstyle{remark}

\def\customauthor{\empty}
\def\customdate{\empty}
\let\oldauthor\author
\renewcommand{\author}[1]{\def\customauthor{#1}}
\renewcommand{\date}[1]{\def\customdate{#1}}
\oldauthor{\customauthor\\\customdate}

\makeatletter
\def\@settitle{\begin{center}%
  \baselineskip10\p@\relax
    \normalfont\LARGE%<- NEW
  \@title
  \end{center}%
}
\makeatother

\makeatletter
\renewcommand{\@secnumfont}{\bfseries}
\makeatother

\begin{document}

\overfullrule=0pt

\author{\bigskip
\bigskip
Dirk Veestraeten\\\\\\}
\date{\small 6 August 2019\\\\}

\address{Amsterdam School of Economics \\
University of Amsterdam\\
Roetersstraat 11\\
1018WB Amsterdam\\
 the Netherlands} 
 
\email{d.j.m.veestraeten@uva.nl}

\title[first passage time density of the Ornstein--Uhlenbeck process \dots ] 
      {\textbf{A new integral equation for the first passage time density of the Ornstein--Uhlenbeck process}}

\begin{abstract}
\noindent The Laplace transform of the first passage time density of the
Ornstein--Uhlenbeck process for a constant threshold contains a ratio of two
parabolic cylinder functions for which no analytical inversion formula is
available. Recently derived inverse Laplace transforms for the product of
two parabolic cylinder functions together with the convolution theorem of
the Laplace transform then allow to derive a new Volterra integral equation
for this first passage time density. The kernel of this integral equation
contains a parabolic cylinder function and the Fortet renewal equation for
the Ornstein--Uhlenbeck process emerges as a special case, namely when
the order $q$ of the parabolic cylinder function is set at $0$. The integral
equation is shown to hold both for constant as well as time dependent
thresholds. Moreover, the kernel of the integral equation is regular for $%
q\leqslant -1$.

\end{abstract}

\maketitle

\noindent  \\

\noindent  \\

\setlength{\baselineskip}{2\baselineskip}

\section{\textbf{Introduction}}
\noindent First passage time probability density functions (pdf) play a prominent role
in various fields such as neurophysiology (see \cite{g67} and the overview
in \cite{sg13}), psychology \cite{sr04}, climate studies \cite{dka05},
biochemistry \cite{kc03} and finance \cite{gy03}. However, closed form
expressions for the first passage time pdf can only be obtained in a very
limited number of cases, see \cite{sg13}.

Hence, applications rely on numerical procedures that are based on a variety of methods.
For instance, \cite{app05} discusses three such methods for the
Ornstein--Uhlenbeck process, namely an eigenfunction expansion, an approach
that relies on the link between the Ornstein--Uhlenbeck process and the
three dimensional Bessel bridge, and an integral representation consisting
of cosine transforms. An alternative and frequently employed integral
representation is the renewal equation of Fortet \cite{f43}. This approach
allows to infer the first passage time pdf when the transition pdf is known
(see\cite{dd08}, \cite{ld08} and \cite{sg13} for more detail). Numerical
applications based on this method are complicated by the fact that the
kernel of this Volterra equation of the first kind is weakly singular, see 
\cite{sg13}. However, this problem can be circumvented, for instance, by
transforming the Fortet renewal equation into a Volterra equation of the
second kind, see \cite{rss84}, \cite{bnr87}, \cite{gnrs89}.

This paper derives a new Volterra integral equation for the first passage
time pdf of the Ornstein--Uhlenbeck process of which the properties may be
attractive for numerical applications. We start from the inverse Laplace
transforms for products of two parabolic cylinder functions that were
recently derived in \cite{v15} and \cite{v17}. As the Laplace transform of
the first passage time pdf for a constant threshold contains a ratio of two
parabolic cylinders (see \cite{rs69}, \cite{smp70}, \cite{cr71}), the
convolution theorem of the Laplace transform then can be used to derive a
new integral equation for this first passage time pdf, namely one in which
the kernel contains a parabolic cylinder function. This integral equation is
a Volterra equation of the first kind, with the exception of a very
specific, but useful, combination of parameters for which it is of the
second kind. The kernel of this integral equation is regular for $q\leqslant
-1$, where $q$ denotes the order of the parabolic cylinder function. The
Fortet renewal equation for the Ornstein--Uhlenbeck process is a special
case of the new integral equation as it emerges for $q=0$. As the Fortet
renewal equation holds both for constant as well as for time depending
thresholds, we verify whether this property also holds for the new integral
equation for values of $q$ that differ from $0$. This indeed is the case and
is verified numerically by using the fact that the first passage time pdf of
the Ornstein--Uhlenbeck process can be obtained in analytical form for
specific exponential time dependent thresholds as noted in \cite{r77} and
implemented in \cite{v14}.

The remainder of this paper is organized as follows. Section $2$ presents
the notation and definitions. Section $3$ uses the convolution theorem of
the Laplace transform to derive an integral equation for the first passage
time pdf of the Ornstein--Uhlenbeck process of which the kernel contains a
parabolic cylinder function of order $q$. For $q=0$, the integral equation
simplifies into the Fortet renewal equation for the Ornstein--Uhlenbeck
process. It is also shown that the integral equation holds for constant as
well as time dependent thresholds and that the kernel is regular for $%
q\leqslant -1$. Section $4$ illustrates two approaches through which the
only known closed form expression for the first passage time pdf of the
Ornstein--Uhlenbeck process under a constant boundary can easily be obtained
from the results in this paper.

\medskip
\section{\textbf{Notation and preliminaries}} 
\noindent Let $\left\{ X\left( t\right) ,t\geqslant 0\right\} $ with $X\left( 0\right)
=x_{0}$ be an Ornstein--Uhlenbeck process with drift $\left( -x/\theta +\mu
\right) $, $\theta >0$, and infinitesimal variance $\sigma ^{2}$. Its
transition pdf is defined as%
\begin{equation*}
f\left( x,t|x_{0}\right) =\dfrac{\partial }{\partial x}P\left\{ X\left(
t\right) \leqslant x|X\left( 0\right) =x_{0}\right\} ,
\end{equation*}%
\noindent which is given by%
\begin{equation}
f\left( x,t|x_{0}\right) =\dfrac{1}{\sqrt{\sigma ^{2}\pi \theta \left(
1-\exp \left( -2t/\theta \right) \right) }}\exp \left( -\dfrac{\left( \left(
x-\mu \theta \right) -\left( x_{0}-\mu \theta \right) \exp \left( -t/\theta
\right) \right) ^{2}}{\sigma ^{2}\theta \left( 1-\exp \left( -2t/\theta
\right) \right) }\right) ,  \label{transition pdf}
\end{equation}%
\noindent see Equation (5.29) in \cite{sg13}.

The first passage time at the constant threshold $S$ with $S>x_{0}$ is the
random variable that is defined by%
\begin{equation*}
T_{S}=\underset{t\geqslant 0}{\inf }\left\{ t:X\left( t\right) >S|X\left(
0\right) =x_{0}\right\} \hspace{0.5cm}x_{0}<S,
\end{equation*}%
\noindent and the first passage time pdf\ is defined as%
\begin{equation*}
g\left( S,t|x_{0}\right) =\dfrac{\partial }{\partial t}P\left\{
T_{S}<t\right\} .
\end{equation*}%
\noindent For the time dependent threshold $S\left( t\right) $ with $S\left(
t\right) >x_{0}$, the following definitions apply%
\begin{equation*}
T_{S\left( t\right) }=\underset{t\geqslant 0}{\inf }\left\{ t:X\left(
t\right) >S\left( t\right) |X\left( 0\right) =x_{0}\right\} \hspace{0.5cm}%
x_{0}<S\left( t\right) ,
\end{equation*}%
\noindent and%
\begin{equation*}
g\left( S\left( t\right) ,t|x_{0}\right) =\dfrac{\partial }{\partial t}%
P\left\{ T_{S\left( t\right) }<t\right\} .
\end{equation*}%
\noindent The Laplace transform of the first passage time pdf for the
constant threshold $S$ with $x_{0}<S$ is denoted by $g_{\lambda }\left(
S\left\vert x_{0}\right. \right) $ and for the Ornstein--Uhlenbeck process
is given by%
\begin{equation}
g_{\lambda }\left( S\left\vert x_{0}\right. \right) =\exp \left( \dfrac{%
\left( x_{0}-\mu \theta \right) ^{2}-\left( S-\mu \theta \right) ^{2}}{%
2\sigma ^{2}\theta }\right) \dfrac{D_{-\lambda \theta }\left( \sqrt{\dfrac{2%
}{\sigma ^{2}\theta }}\left( \mu \theta -x_{0}\right) \right) }{D_{-\lambda
\theta }\left( \sqrt{\dfrac{2}{\sigma ^{2}\theta }}\left( \mu \theta
-S\right) \right) }\hspace{0.5cm}x_{0}<S,  \label{lt fpt}
\end{equation}%
\noindent see Equation (5.37) in \cite{sg13}, where $D_{\nu }\left( z\right) 
$ denotes the parabolic cylinder function of order $\nu $ and argument $z$,
see \cite{oms09}.

The Fortet renewal equation \cite{f43} is a Volterra integral equation of
the first kind that connects the transition pdf and the first passage time
pdf for constant $S$ as follows%
\begin{equation}
f\left( x,t|x_{0}\right) =\int_{0}^{t}f\left( x,t|S,\tau \right) g\left(
S,\tau |x_{0}\right) d\tau \hspace{0.5cm}x\geqslant S,  \label{Fortet S}
\end{equation}%
\noindent and for time dependent $S\left( t\right) $ via 
\begin{equation}
f\left( x,t|x_{0}\right) =\int_{0}^{t}f\left( x,t|S\left( \tau \right) ,\tau
\right) g\left( S\left( \tau \right) ,\tau |x_{0}\right) d\tau \hspace{0.5cm}%
x\geqslant S\left( t\right) ,  \label{Fortet S(t)}
\end{equation}%
\noindent see \cite{sg13}. Specializing the Fortet equations (\ref{Fortet S}%
) and (\ref{Fortet S(t)}) for the Ornstein--Uhlenbeck process, i.e. when
using the transition density (\ref{transition pdf}), gives%
\begin{eqnarray}
&&\int\nolimits_{0}^{t}\dfrac{1}{\sqrt{\sigma ^{2}\pi \theta \left( 1-\exp
\left( -2\left( t-\tau \right) /\theta \right) \right) }}\exp \left( -\dfrac{%
\left( \left( x-\mu \theta \right) -\left( S-\mu \theta \right) \exp \left(
-\left( t-\tau \right) /\theta \right) \right) ^{2}}{\sigma ^{2}\theta
\left( 1-\exp \left( -2\left( t-\tau \right) /\theta \right) \right) }%
\right)   \notag \\
&&\hspace{1cm}\times g\left( S,\tau \left\vert x_{0}\right. \right) d\tau =
\label{Fortet S OU} \\
&&\hspace{1cm}\dfrac{1}{\sqrt{\sigma ^{2}\pi \theta \left( 1-\exp \left(
-2t/\theta \right) \right) }}\exp \left( -\dfrac{\left( \left( x-\mu \theta
\right) -\left( x_{0}-\mu \theta \right) \exp \left( -t/\theta \right)
\right) ^{2}}{\sigma ^{2}\theta \left( 1-\exp \left( -2t/\theta \right)
\right) }\right) \hspace{0.5cm}x\geqslant S,  \notag
\end{eqnarray}%
\noindent and%
\begin{eqnarray}
&&\int\nolimits_{0}^{t}\frac{1}{\sqrt{\sigma ^{2}\pi \theta \left( 1-\exp
\left( -2\left( t-\tau \right) /\theta \right) \right) }}\exp \left( -\dfrac{%
\left( \left( x-\mu \theta \right) -\left( S\left( \tau \right) -\mu \theta
\right) \exp \left( -\left( t-\tau \right) /\theta \right) \right) ^{2}}{%
\sigma ^{2}\theta \left( 1-\exp \left( -2\left( t-\tau \right) /\theta
\right) \right) }\right)   \notag \\
&&\hspace{1cm}\times g\left( S\left( \tau \right) ,\tau \left\vert
x_{0}\right. \right) d\tau =  \label{Fortet S(t) OU} \\
&&\hspace{1cm}\dfrac{1}{\sqrt{\sigma ^{2}\pi \theta \left( 1-\exp \left(
-2t/\theta \right) \right) }}\exp \left( -\dfrac{\left( \left( x-\mu \theta
\right) -\left( x_{0}-\mu \theta \right) \exp \left( -t/\theta \right)
\right) ^{2}}{\sigma ^{2}\theta \left( 1-\exp \left( -2t/\theta \right)
\right) }\right) \hspace{0.5cm}x\geqslant S\left( t\right) .  \notag
\end{eqnarray}%
\noindent The Laplace transforms of the original functions $f_{1}\left(
t\right) $ and $f_{2}\left( t\right) $ are denoted by%
\begin{eqnarray*}
\overline{f}_{1}\left( \lambda \right)  &=&\int_{0}^{\infty }\exp \left(
-\lambda t\right) f_{1}\left( t\right) dt, \\
\overline{f}_{2}\left( \lambda \right)  &=&\int_{0}^{\infty }\exp \left(
-\lambda t\right) f_{2}\left( t\right) dt,
\end{eqnarray*}%
\noindent where $\operatorname{Re}\lambda >0$. The convolution theorem of the
Laplace transform then states%
\begin{equation}
\overline{f}_{1}\left( \lambda \right) \overline{f}_{2}\left( \lambda
\right) =\int_{0}^{\infty }\exp \left( -\lambda t\right) f_{1}\left(
t\right) \ast f_{2}\left( t\right) dt,  \label{lt convol1}
\end{equation}%
\noindent where $f_{1}\left( t\right) \ast f_{2}\left( t\right) $ is the
convolution of $f_{1}\left( t\right) $ and $f_{2}\left( t\right) $ that is
to be obtained from%
\begin{equation}
f_{1}\left( t\right) \ast f_{2}\left( t\right) =\int_{0}^{t}f_{1}\left( \tau
\right) f_{2}\left( t-\tau \right) d\tau ,  \label{ilt convol1}
\end{equation}%
\noindent see \cite{db15}.

\medskip
\section{\textbf{The integral equation and its properties}}
\noindent The first passage time pdf $g\left( S,t\left\vert x_{0}\right.
\right) $ and its Laplace transform $g_{\lambda }\left( S\left\vert
x_{0}\right. \right) $ will represent $f_{1}\left( t\right) $ and $\overline{%
f}_{1}\left( \lambda \right) $, respectively, within the convolution. Thus, $%
\overline{f}_{1}\left( \lambda \right) $ contains a ratio of parabolic
cylinder functions, see (\ref{lt fpt}). The initial function $f_{2}\left(
t\right) $ then will need to be chosen in such manner that the inversion
formula for $\overline{f}_{1}\left( \lambda \right) \overline{f}_{2}\left(
\lambda \right) $ is known.

The paper starts from the inverse Laplace transform for the product of two
parabolic cylinder functions with different arguments and orders that was
recently derived in \cite{v17} and that generalized the results in \cite{v15}%
. The inverse Laplace transform in Equation (2.1) in \cite{v17} first will
be extended via the result in Equation (4.7) in \cite{v15}. This gives the
following inverse Laplace transform in which attention is restricted to real
arguments and orders of the parabolic cylinder functions%
\begin{eqnarray}
&&\Gamma \left[ \lambda +c\right] D_{-c-\lambda }\left( y\right)
D_{q-\lambda }\left( z\right) -\left. \sqrt{\dfrac{\pi }{2}}\right\vert
_{y+z=0,c+q=1}=  \notag \\
&&\int\nolimits_{0}^{\infty }\exp \left( -\left( \lambda +c\right) t\right)
\left( 1-\exp \left( -2t\right) \right) ^{-\left( 1+c+q\right) /2}\exp
\left( -\dfrac{\left( y+z\exp \left( -t\right) \right) ^{2}}{4\left( 1-\exp
\left( -2t\right) \right) }\right)   \notag \\
&&\hspace{1cm}\times D_{c+q}\left( \dfrac{z+y\exp \left( -t\right) }{\sqrt{%
1-\exp \left( -2t\right) }}\right) dt  \label{ilt double1} \\[0.3cm]
&&\left[ \operatorname{Re}\left( \lambda +c\right) >0\text{ and }y+z>0\text{ or }%
y+z=0,c+q\leqslant 1\right] .  \notag
\end{eqnarray}%

\noindent The term $\left. \sqrt{\pi /2}\right\vert _{_{y+z=0,c+q=1}}$ in
the inverse Laplace transform (\ref{ilt double1}) equals $\sqrt{\pi /2}$
when $y+z=0$ and $c+q=1$, and equals $0$ otherwise. This additional term
finds its origin in the use of the differentiation property of the Laplace
transform in the derivation of the inverse Laplace transform in Equation
(4.7) in \cite{v15}. This extra term turns out to be quite useful as it
later allows to obtain the only known closed form expression for the first
passage time pdf under a constant threshold in a very simple manner, see
Section $4$. For $y+z=0$ and $c+q>1$, the term $f\left( 0\right) $ within
the application of the differentiation property of the Laplace transform in 
\cite{v15} becomes infinite such that the inverse Laplace transform (\ref%
{ilt double1}) requires the restriction $c+q\leqslant 1$ when $y+z=0$, but
no such restriction on $c+q$ is to be imposed when $y+z>0$.

The time scaling property of the Laplace transform states 
\begin{equation}
\overline{f}\left( \lambda \theta \right) =\frac{1}{\theta }L\left\{ f\left( 
\frac{t}{\theta }\right) \right\} \hspace{0.5cm}\theta >0,
\label{time scaling}
\end{equation}%
\noindent see Equation (29.2.13) in \cite{as72}. Using the latter property, $%
c=0$, $y=\sqrt{\dfrac{2}{\sigma ^{2}\theta }}\left( \mu \theta -S\right) $
and $z=\sqrt{\dfrac{2}{\sigma ^{2}\theta }}\left( x-\mu \theta \right) $
allows to rewrite (\ref{ilt double1}) as follows%
\begin{eqnarray}
&&\Gamma \left[ \lambda \theta \right] D_{-\lambda \theta }\left( \sqrt{%
\dfrac{2}{\sigma ^{2}\theta }}\left( \mu \theta -S\right) \right)
D_{q-\lambda \theta }\left( \sqrt{\dfrac{2}{\sigma ^{2}\theta }}\left( x-\mu
\theta \right) \right) -\left. \sqrt{\dfrac{\pi }{2}}\right\vert _{x=S,q=1}=
\notag \\
&&\int\nolimits_{0}^{\infty }\exp \left( -\lambda t\right) \left\{ \dfrac{1%
}{\theta }\left( 1-\exp \left( -2t/\theta \right) \right) ^{-\left(
1+q\right) /2}\exp \left( -\dfrac{\left( \left( \mu \theta -S\right) +\left(
x-\mu \theta \right) \exp \left( -t/\theta \right) \right) ^{2}}{2\sigma
^{2}\theta \left( 1-\exp \left( -2t/\theta \right) \right) }\right) \right. 
\notag \\
&&\hspace{1cm}\left. \times D_{q}\left( \dfrac{\sqrt{2}\left( \left( x-\mu
\theta \right) +\left( \mu \theta -S\right) \exp \left( -t/\theta \right)
\right) }{\sqrt{\sigma ^{2}\theta \left( 1-\exp \left( -2t/\theta \right)
\right) }}\right) \right\} dt  \label{ilt double2} \\[0.3cm]
&&\left[ \operatorname{Re}\lambda >0\text{ and }x>S\text{ or }x=S,q\leqslant 1\right]
.  \notag
\end{eqnarray}%

\noindent The two terms on the left hand side of the Laplace transform (\ref%
{ilt double2}) together specify the Laplace transform $\overline{f}%
_{2}\left( \lambda \right) $ that will be used within the convolution. The
term within curly brackets on the right hand side of the Laplace transform (%
\ref{ilt double2}) then represents the original function $f_{2}\left(
t\right) $.

The product $\overline{f}_{1}\left( \lambda \right) \overline{f}_{2}\left(
\lambda \right) $, given the above choices for $\overline{f}_{1}\left(
\lambda \right) $ and $\overline{f}_{2}\left( \lambda \right) $, then
emerges as%
\begin{eqnarray}
&&\overline{f}_{1}\left( \lambda \right) \overline{f}_{2}\left( \lambda
\right) =\Gamma \left[ \lambda \theta \right] \exp \left( \dfrac{\left(
x_{0}-\mu \theta \right) ^{2}-\left( S-\mu \theta \right) ^{2}}{2\sigma
^{2}\theta }\right) D_{-\lambda \theta }\left( \sqrt{\dfrac{2}{\sigma
^{2}\theta }}\left( \mu \theta -x_{0}\right) \right)   \notag \\
&&\hspace{1cm}\times D_{q-\lambda \theta }\left( \sqrt{\dfrac{2}{\sigma
^{2}\theta }}\left( x-\mu \theta \right) \right) -\left. \sqrt{\dfrac{\pi }{2%
}}\right\vert _{x=S,q=1}g_{\lambda }\left( S\left\vert x_{0}\right. \right) 
\hspace{0.5cm}x_{0}<S.  \label{lt convol2}
\end{eqnarray}%

\noindent Employing the inverse Laplace transform (\ref{ilt double1}) and
using the above values for $c$, $y$ and $z$ allow to express the inversion
formula for $\overline{f}_{1}\left( \lambda \right) \overline{f}_{2}\left(
\lambda \right) $ as%
\begin{eqnarray}
&&\int\nolimits_{0}^{\infty }\exp \left( -\lambda t\right) \left\{ \exp
\left( \dfrac{\left( x_{0}-\mu \theta \right) ^{2}-\left( S-\mu \theta
\right) ^{2}}{2\sigma ^{2}\theta }\right) \dfrac{1}{\theta }\left( 1-\exp
\left( -2t/\theta \right) \right) ^{-\left( 1+q\right) /2}\right.   \notag \\
&&\hspace{1cm}\times \exp \left( -\dfrac{\left( \left( \mu \theta
-x_{0}\right) +\left( x-\mu \theta \right) \exp \left( -t/\theta \right)
\right) ^{2}}{2\sigma ^{2}\theta \left( 1-\exp \left( -2t/\theta \right)
\right) }\right)   \notag \\
&&\hspace{1cm}\left. \times D_{q}\left( \dfrac{\sqrt{2}\left( \left( x-\mu
\theta \right) +\left( \mu \theta -x_{0}\right) \exp \left( -t/\theta
\right) \right) }{\sqrt{\sigma ^{2}\theta \left( 1-\exp \left( -2t/\theta
\right) \right) }}\right) \right\} dt  \notag \\
&&\hspace{1cm}-\left\{ \left. \sqrt{\dfrac{\pi }{2}}\right\vert
_{x=S,q=1}g\left( S,t\left\vert x_{0}\right. \right) \right\} 
\label{ilt convol2} \\[0.3cm]
&&\hspace{0.5cm}\left[ \operatorname{Re}\lambda >0\text{ and }x>S\text{ or }%
x=S,q\leqslant 1\right] .  \notag
\end{eqnarray}%
\noindent Hence, $f_{1}\left( t\right) \ast f_{2}\left( t\right) $ is to be
obtained from the two terms within curly brackets in the inversion formula (%
\ref{ilt convol2}). The convolution theorem (\ref{ilt convol1}) then gives
-- after some straightforward simplifications -- the following integral
equation for the first passage time pdf of the Ornstein--Uhlenbeck process
for the constant threshold $S$%
\begin{eqnarray}
&&\int\nolimits_{0}^{t}\left( 1-\exp \left( -2\left( t-\tau \right) /\theta
\right) \right) ^{-\left( 1+q\right) /2}\exp \left( -\dfrac{\left( \left(
x-\mu \theta \right) -\left( S-\mu \theta \right) \exp \left( -\left( t-\tau
\right) /\theta \right) \right) ^{2}}{2\sigma ^{2}\theta \left( 1-\exp
\left( -2\left( t-\tau \right) /\theta \right) \right) }\right)   \notag \\
&&\hspace{1cm}\times D_{q}\left( \dfrac{\sqrt{2}\left( \left( x-\mu \theta
\right) -\left( S-\mu \theta \right) \exp \left( -\left( t-\tau \right)
/\theta \right) \right) }{\sqrt{\sigma ^{2}\theta \left( 1-\exp \left(
-2\left( t-\tau \right) /\theta \right) \right) }}\right) g\left( S,\tau
\left\vert x_{0}\right. \right) d\tau =  \notag \\
&&\left( 1-\exp \left( -2t/\theta \right) \right) ^{-\left( 1+q\right)
/2}\exp \left( -\dfrac{\left( \left( x-\mu \theta \right) -\left( x_{0}-\mu
\theta \right) \exp \left( -t/\theta \right) \right) ^{2}}{2\sigma
^{2}\theta \left( 1-\exp \left( -2t/\theta \right) \right) }\right) 
\label{integral equation1} \\
&&\hspace{1cm}\times D_{q}\left( \dfrac{\sqrt{2}\left( \left( x-\mu \theta
\right) -\left( x_{0}-\mu \theta \right) \exp \left( -t/\theta \right)
\right) }{\sqrt{\sigma ^{2}\theta \left( 1-\exp \left( -2t/\theta \right)
\right) }}\right) -\theta \left. \sqrt{\dfrac{\pi }{2}}\right\vert
_{x=S,q=1}g\left( S,t\left\vert x_{0}\right. \right)   \notag \\[0.3cm]
&&\left[ x>S\text{ or }x=S,q\leqslant 1\right] .  \notag
\end{eqnarray}%

\noindent The kernel of the Volterra integral equation (\ref{integral
equation1}) is a combination of terms that contain exponential functions as
well as a parabolic cylinder function of order $q$. Note that no restriction
on $q$ is required for $x>S$, whereas $q\leqslant 1$ is to be imposed when $%
x=S$. This integral equation, notwithstanding the fact that it relies on an
argument for constant $S$, actually also holds for the time dependent
threshold $S\left( t\right) $. In order to illustrate this, we first
specialize the integral equation (\ref{integral equation1}) in terms of the
Fortet equation (\ref{Fortet S OU}). Using $q=0$ simplifies the parabolic
cylinder function into an exponential function via the following identity%
\begin{equation*}
D_{0}\left( z\right) =\exp \left( -\dfrac{z^{2}}{4}\right) ,
\end{equation*}%
\noindent see Equation (46:4:1) in \cite{oms09}. Plugging the latter
property into the integral equation (\ref{integral equation1}) and
simplifying then indeed gives the Fortet equation for the
Ornstein--Uhlenbeck process with constant threshold $S$ that is presented in
expression (\ref{Fortet S OU}). The renewal equation (\ref{Fortet S OU})
extends also to the time dependent threshold $S\left( t\right) $ resulting
in the renewal equation (\ref{Fortet S(t) OU}). Or, the specialization of
the integral equation (\ref{integral equation1}) for $q=0$ likewise extends
to the time dependent threshold $S\left( t\right) $. The integral equations (%
\ref{Fortet S OU}) and (\ref{Fortet S(t) OU}) then raise the question
whether or not this extension of applicability from constant to time
dependent threshold also holds for values of $q$ that differ from $0$. The
answer is affirmative and this will be illustrated for the case of a time
dependent threshold for which the first passage time pdf is known in closed
form. As noted in \cite{r77}, a closed form expression for the first passage
time pdf of the Ornstein--Uhlenbeck process can be obtained for thresholds
that depend on time in a specific exponential manner. In fact, \cite{r77}
uses the Doob transformation \cite{d42} to express the problem for the
Ornstein--Uhlenbeck process in terms of the Wiener process with a constant
or affine threshold for which closed form expressions exist. The resulting
transition and first passage time pdf are illustrated in \cite{v14} for the
threshold $S\left( t\right) _{\exp }$ with%
\begin{equation}
S\left( t\right) _{\exp }=d_{1}\exp \left( -t/\theta \right) +d_{2}\sinh
\left( t/\theta \right) +\mu \theta \left( 1-\exp \left( -t/\theta \right)
\right) ,  \label{threshold exp}
\end{equation}%
\noindent see Equation (2.2) in \cite{v14}, where $d_{1}$ and $d_{2}$ are
arbitrary constants. The first passage time pdf for the latter time
dependent threshold is%
\begin{eqnarray}
&&g\left( S\left( t\right) _{\exp },t\left\vert x_{0}\right. \right) =\dfrac{%
2\left( d_{1}-x_{0}\right) }{\sigma \theta ^{3/2}\sqrt{\pi }}\dfrac{\exp
\left( -t/\theta \right) }{\left( 1-\exp \left( -2t/\theta \right) \right)
^{3/2}}  \notag \\
&&\hspace{1cm}\times \exp \left( -\dfrac{\left( \left( d_{1}-x_{0}\right)
\exp \left( -t/\theta \right) +d_{2}\sinh \left( t/\theta \right) \right)
^{2}}{\sigma ^{2}\theta \left( 1-\exp \left( -2t/\theta \right) \right) }%
\right) \hspace{0.5cm}x_{0}<S\left( t\right) _{\exp },  \label{fpt exp}
\end{eqnarray}%
\noindent see Equation (2.13) in \cite{v14}. The applicability of the
integral equation (\ref{integral equation1}) for both constant as well as
time dependent thresholds then can be tested by inserting the time dependent
threshold (\ref{threshold exp}) and its first passage time pdf (\ref{fpt exp}%
) into the integral equation (\ref{integral equation1}). The resulting
integral cannot be evaluated analytically. However, numerical evaluation
reveals that the integral equals the term on the right hand side for all
values of $q$. Hence, the integral equation (\ref{integral equation1}) can
indeed be generalized towards the time dependent threshold $S\left( t\right) 
$. This gives the following integral equation for the first passage time pdf
of the Ornstein--Uhlenbeck process with time dependent threshold $S\left(
t\right) $%
\begin{eqnarray}
&&\int\nolimits_{0}^{t}\left( 1-\exp \left( -2\left( t-\tau \right) /\theta
\right) \right) ^{-\left( 1+q\right) /2}\exp \left( -\dfrac{\left( \left(
x-\mu \theta \right) -\left( S\left( \tau \right) -\mu \theta \right) \exp
\left( -\left( t-\tau \right) /\theta \right) \right) ^{2}}{2\sigma
^{2}\theta \left( 1-\exp \left( -2\left( t-\tau \right) /\theta \right)
\right) }\right)   \notag \\
&&\hspace{1cm}\times D_{q}\left( \dfrac{\sqrt{2}\left( \left( x-\mu \theta
\right) -\left( S\left( \tau \right) -\mu \theta \right) \exp \left( -\left(
t-\tau \right) /\theta \right) \right) }{\sqrt{\sigma ^{2}\theta \left(
1-\exp \left( -2\left( t-\tau \right) /\theta \right) \right) }}\right)
g\left( S\left( \tau \right) ,\tau \left\vert x_{0}\right. \right) d\tau = 
\notag \\
&&\left( 1-\exp \left( -2t/\theta \right) \right) ^{-\left( 1+q\right)
/2}\exp \left( -\dfrac{\left( \left( x-\mu \theta \right) -\left( x_{0}-\mu
\theta \right) \exp \left( -t/\theta \right) \right) ^{2}}{2\sigma
^{2}\theta \left( 1-\exp \left( -2t/\theta \right) \right) }\right) 
\label{integral equation2} \\
&&\hspace{1cm}\times D_{q}\left( \dfrac{\sqrt{2}\left( \left( x-\mu \theta
\right) -\left( x_{0}-\mu \theta \right) \exp \left( -t/\theta \right)
\right) }{\sqrt{\sigma ^{2}\theta \left( 1-\exp \left( -2t/\theta \right)
\right) }}\right) -\theta \left. \sqrt{\dfrac{\pi }{2}}\right\vert
_{x=S\left( t\right) ,q=1}g\left( S\left( t\right) ,t\left\vert x_{0}\right.
\right)   \notag \\[0.3cm]
&&\left[ x>S\left( t\right) \text{ or }x=S\left( t\right) ,q\leqslant 1%
\right] .  \notag
\end{eqnarray}%

\noindent No restrictions apply to the order $q$ in the integral equation (%
\ref{integral equation2}) except when $x=S$ $\left( t\right) $ in which case
the restriction $q\leqslant 1$ is to be imposed. The integral equation (\ref%
{integral equation2}) is a Volterra equation of the first kind for $%
x>S\left( t\right) $ or $q\neq 1$, and is of the second kind when $x=S\left(
t\right) $ and $q=1$.

The kernel in the integral equation (\ref{integral equation2}) is regular
for a wide range of values of the order $q$. First, we examine the behaviour
of the kernel when evaluating the limit $\tau \rightarrow t$ for $x>S\left(
t\right) $. The kernel then goes to zero for all values of $q$. Indeed, the
parabolic cylinder function moves to $0$ given $\underset{z\rightarrow
+\infty }{\lim }\left[ D_{\nu }\left( z\right) \right] =0\,$, see Equation
(46:7:1) in \cite{oms09}. Also, the limit of the second term on the left
hand side of (\ref{integral equation2}) vanishes for $\tau \rightarrow t$.
The first term grows larger than $1$ for $q>-1$, but as the second and third
term go to $0$ at a faster rate, the kernel vanishes for $\tau \rightarrow t$
and $x>S\left( t\right) $ irrespective of the level of $q$. Second, the
argument in the second term goes to $0$ for $\tau \rightarrow t$ when $%
x=S\left( t\right) $. Also the argument in the parabolic cylinder function
then equals $0$ which gives%
\begin{equation}
D_{q}\left( 0\right) =\dfrac{2^{q/2}\sqrt{\pi }}{\Gamma \left( \dfrac{1-q}{2}%
\right) },  \label{D arg 0}
\end{equation}%
\noindent see Equation (46:7:1) in \cite{oms09}. The latter expression is
finite except for $q=1,3,5,7,...$ However, the first term in the kernel
becomes infinite for $q>-1$. Taken together, the kernel\ of the integral
equation (\ref{integral equation2}) is regular for $q\leqslant -1$ and $%
x=S\left( t\right) $ and is weakly singular for $q>-1$ and $x=S\left(
t\right) .$

Finally, it is to be noted that the kernel of the integral equation (\ref%
{integral equation2}) can be expressed in terms of other special functions
and polynomials for specific values of $q$. For positive integer values of $q
$, the parabolic cylinder function simplifies into the Hermite polynomials
via the relation%
\begin{equation*}
D_{n}\left( z\right) =2^{-n/2}\exp \left( -z^{2}/4\right) H_{n}\left( z/%
\sqrt{2}\right) \hspace{0.5cm}n=0,1,2,3,...
\end{equation*}%
\noindent see Equation (46:4:1) in \cite{oms09}. The parabolic cylinder
function reduces into the complementary error function erfc$\left( z\right) $
for negative integers via%
\begin{equation*}
D_{-n-1}\left( z\right) =2^{-1/2}\sqrt{\pi }\dfrac{\left( -1\right) ^{n}}{n!}%
\exp \left( -z^{2}/4\right) \dfrac{d\hspace{0.05cm}^{n}}{d\hspace{0.01cm}%
z^{n}}\left[ \exp \left( z^{2}/2\right) \text{erfc}\left( 2^{-1/2}z\right) %
\right] \hspace{0.5cm}n=0,1,2,3,...
\end{equation*}%
\noindent see p. 326 in \cite{mos66}. The kernel can also be specified in
terms of the modified Bessel function of the second kind given the property%
\begin{equation*}
D_{-1/2}\left( z\right) =\sqrt{z/2\pi }\hspace{0.1cm}K_{1/4}\left(
z^{2}/4\right) ,
\end{equation*}%
\noindent see Equation (46:4:4) in \cite{oms09}.

\medskip
\section{\textbf{Retrieving the first passage time pdf for $\pmb{S=\protect\mu \protect%
\theta} $}}
\noindent The first passage time pdf for the Ornstein--Uhlenbeck process with a
constant threshold at $S=\mu \theta $ possesses the following well known
closed form expression%
\begin{equation}
g\left( \mu \theta ,t\left\vert x_{0}\right. \right) =\dfrac{2\left( \mu
\theta -x_{0}\right) }{\sigma \theta ^{3/2}\sqrt{\pi }}\dfrac{\exp \left(
-t/\theta \right) }{\left( 1-\exp \left( -2t/\theta \right) \right) ^{3/2}}%
\exp \left( -\dfrac{\left( x_{0}-\mu \theta \right) ^{2}\exp \left(
-2t/\theta \right) }{\sigma ^{2}\theta \left( 1-\exp \left( -2t/\theta
\right) \right) }\right) \hspace{0.5cm}x_{0}<\mu \theta ,  \label{S=mt}
\end{equation}%
\noindent see \cite{smp70} and \cite{r77}, who obtained this result via
Doob's transformation \cite{d42}. However, the results from the previous
section can be used to obtain the first passage time pdf (\ref{S=mt}) in two
ways that do not rely on the latter transformation.

First, the integral equation (\ref{integral equation2}) simplifies into the
following Volterra equation of the second kind for $x=S$ and $q=1$%
\begin{eqnarray}
&&g\left( x,t\left\vert x_{0}\right. \right) =  \notag \\
&&\dfrac{2\left( \left( x-\mu \theta \right) -\left( x_{0}-\mu \theta
\right) \exp \left( -t/\theta \right) \right) }{\sigma \theta ^{3/2}\sqrt{%
\pi }\left( 1-\exp \left( -2t/\theta \right) \right) ^{3/2}}\exp \left( -%
\dfrac{\left( \left( x-\mu \theta \right) -\left( x_{0}-\mu \theta \right)
\exp \left( -t/\theta \right) \right) ^{2}}{\sigma ^{2}\theta \left( 1-\exp
\left( -2t/\theta \right) \right) }\right)   \notag \\
&&-\int\nolimits_{0}^{t}\dfrac{2\left( x-\mu \theta \right) \left( 1-\exp
\left( -\left( t-\tau \right) /\theta \right) \right) }{\sigma \theta ^{3/2}%
\sqrt{\pi }\left( 1-\exp \left( -2\left( t-\tau \right) /\theta \right)
\right) ^{3/2}}  \label{integral equation3} \\
&&\hspace{1cm}\times \exp \left( -\dfrac{\left( x-\mu \theta \right)
^{2}\left( 1-\exp \left( -\left( t-\tau \right) /\theta \right) \right) ^{2}%
}{\sigma ^{2}\theta \left( 1-\exp \left( -2\left( t-\tau \right) /\theta
\right) \right) }\right) g\left( x,\tau \left\vert x_{0}\right. \right)
d\tau \hspace{1cm}x_{0}<x  \notag
\end{eqnarray}%
\noindent of which the derivation used the property 
\begin{equation}
D_{1}\left( z\right) =z\exp \left( -\dfrac{z^{2}}{4}\right) ,  \label{D1}
\end{equation}%
\noindent see Equation (46:4:1) in \cite{oms09}. Using $x=\mu \theta $, i.e. 
$S=\mu \theta $, causes the integral in the integral equation (\ref{integral
equation3}) to vanish such that indeed the first passage time pdf (\ref{S=mt}%
) emerges.

Second, plugging $S=\mu \theta $ into the Laplace transform of the first
passage time pdf (\ref{lt fpt}) and using the property (\ref{D arg 0}) gives%
\begin{equation}
g_{\lambda }\left( \mu \theta \left\vert x_{0}\right. \right) =\dfrac{1}{%
\sqrt{\pi }}\exp \left( \dfrac{\left( x_{0}-\mu \theta \right) ^{2}}{2\sigma
^{2}\theta }\right) 2^{\lambda \theta /2}\Gamma \left[ \dfrac{1+\lambda
\theta }{2}\right]
D_{-\lambda \theta }\left( \sqrt{\dfrac{2}{\sigma
^{2}\theta }}\left( \mu \theta -x_{0}\right) \right) \hspace{0.5cm}x_{0}<\mu
\theta .  \label{lt fpt S=mt}
\end{equation}%
\noindent The inversion formula for the Laplace transform (\ref{lt fpt S=mt}%
) can be obtained from the inverse Laplace transform (\ref{ilt double1})
when using $y=0$, applying the time scaling property (\ref{time scaling})
and setting $c$ and $q$ at $1$ and $0$, respectively. Subsequently, the
properties (\ref{D arg 0}) and (\ref{D1}) are to be used and further
simplifications rely on the Legendre duplication formula for the gamma
function%
\begin{equation*}
\Gamma \left[ 2z\right] =\left( 2\pi \right) ^{-1/2}2^{2z-1/2}\Gamma \left[ z%
\right] \Gamma \left[ 1/2+z\right] ,
\end{equation*}%
\noindent see Equation (6.1.18) in \cite{as72}. These simplifications then
yield the following inverse Laplace transform%
\begin{eqnarray}
&&2^{\theta \lambda /2}\Gamma \left[ \dfrac{1+\lambda \theta }{2}\right]
D_{-\lambda \theta }\left( z\right) -\left. \sqrt{\pi }\right\vert
_{z=0}=\int\nolimits_{0}^{\infty }\exp \left( -\lambda t\right) \dfrac{%
\sqrt{2}}{\theta }z\exp \left( -t/\theta \right)   \notag \\
&&\hspace{1cm}\times \left( 1-\exp \left( -2t/\theta \right) \right)
^{-3/2}\exp \left( -\dfrac{z^{2}\left( 1+\exp \left( -2t/\theta \right)
\right) }{4\left( 1-\exp \left( -2t/\theta \right) \right) }\right) dt
\label{ilt single} \\
&&\left[ \operatorname{Re}\left( \lambda +1/\theta \right) >0,z\geqslant 0\right] . 
\notag
\end{eqnarray}%
\noindent Specializing the inverse Laplace transform (\ref{ilt single}) for $%
z=\sqrt{\dfrac{2}{\sigma ^{2}\theta }}\left( \mu \theta -x_{0}\right) ,$
i.e. for $z>0$, and using the resulting inversion formula for the Laplace
transform (\ref{lt fpt S=mt}) then gives the first passage time pdf (\ref%
{S=mt}).

\end{document}